\documentclass[12pt]{amsart}
\usepackage[utf8]{inputenc}

\usepackage{subfiles}
\usepackage{hyperref}
\usepackage{cleveref}
\usepackage{amsmath}
\usepackage{amssymb}
\usepackage{amsthm}
\usepackage{mathtools}
\usepackage{thmtools}
\usepackage{thm-restate}
\usepackage[margin=1in]{geometry}
\usepackage{cleveref}
\usepackage{ytableau}
\usepackage{comment}
\usepackage{soul}
\usepackage{tikz}
\usepackage{tikzscale}
\usetikzlibrary{positioning}
\usetikzlibrary{patterns}
\usepackage{graphicx}
\usepackage{caption}
\usepackage{subcaption}
\usepackage{float}

\usepackage[numbers]{natbib}
\numberwithin{equation}{section}

\theoremstyle{definition}
\newtheorem{theorem}{Theorem}[section]
\newtheorem*{theorem*}{Theorem}

\newtheorem*{example*}{Example}
\newtheorem{lemma}[theorem]{Lemma}
\newtheorem*{lemma*}{Lemma}
\newtheorem{corollary}[theorem]{Corollary}
\newtheorem*{corollary*}{Corollary}

\newtheorem*{definition*}{Definition}
\newtheorem{proposition}[theorem]{Proposition}
\newtheorem*{proposition*}{Proposition}

\newtheorem*{remark*}{Remark}
\newtheorem{conjecture}[theorem]{Conjecture}

\renewcommand{\S}{\mathcal S}
\newcommand{\s}{\sigma}
\newcommand{\av}{\text{Av}}
\newcommand{\sort}{\text{Sort}}
\renewcommand{\sc}{\text{SC}}

\newcommand{\comp}{\text{comp}}
\newcommand{\rev}{\text{rev}}
\newcommand{\un}[1]{\underline{#1}}
\newcommand{\sy}{\mathfrak S}
\usepackage{ mathrsfs }

\title{Stack-sorting with stacks avoiding vincular patterns}

\author{William Zhao}\address{\textsc{Dougherty Valley High School}, 10550 Albion Road, San Ramon, CA, USA} \email{williamzsrv@gmail.com}

\begin{document}















\begin{abstract}
We introduce the stack-sorting map $\sc_\s$ that sorts, in a right-greedy manner, an input permutation through a stack that avoids some vincular pattern $\s$. The stack-sorting maps of Cerbai et al. in which the stack avoids a pattern classically and Defant and Zheng in which the stack avoids a pattern consecutively follow as special cases. We first characterize and enumerate the sorting class $\sort(\sc_\s)$, the set of permutations sorted by $s\circ\sc_\s$, for seven length $3$ patterns $\s$. We also decide when $\sort(\sc_\s)$ is a permutation class. Next, we compute $\max_{\pi\in \sy_n}|\sc_\s^{-1}(\pi)|$ and characterize the periodic points of $\sc_\s$ for several length $3$ patterns $\s$. We end with several conjectures and open problems.

\end{abstract}


\maketitle


\footnotesize


\normalsize

\section{Introduction} \label{intro}
Let $\sy_n$ be the symmetric group. Consider \citet{West}'s stack-sorting map $s:\sy_n\to \sy_n$. A stack, which is required to be increasing, processes the input permutation $\tau\in \sy_n$ in a right-greedy manner. In other words, the map $s$ uses a stack that avoids the pattern $21$. In 2019, \citet{CCF} generalized West's stack-sorting map. The map $s_\s:\sy_n\to \sy_n$ processes the input permutation through a stack in a similar manner as $s$, but with a stack that must avoid some pattern $\s$, as in \Cref{exs123}. Since then, the map $s_\s$ and variants have been studied extensively \cite{KatalinBerlow, CCF, 132machine, Defant}. Most relevant to our work is \citet{Defant} that introduced the consecutive-pattern-avoiding stack-sorting map $\sc_\s:\sy_n\to \sy_n$, a modification of $s_\s$ where the stack must consecutively avoid $\s$, as in \Cref{exsc123}. 

\medskip

\citet{BabsonSteinVincular} first introduced vincular pattern avoidance. When considering whether a permutation $\pi$ contains a vincular pattern $\s$, some elements may be required to be adjacent in $\pi$, as indicated by underlined terms in $\s$. For instance, the pattern $1423$ contains $1\un{23}$ and $123$, but avoids $\un{12}3$ and $\un{123}$. Note that classical and consecutive pattern avoidance are both special cases of vincular pattern avoidance. 

\medskip

Therefore, a natural extension of \citet{CCF} and \citet{Defant}'s maps is the generalized map $\sc_\s:\sy_n\to \sy_n$ that sorts with a stack which avoids a vincular pattern $\s$. Consider the permutation $\tau=514362$. \Cref{exs123,exsc123,exscs123,exssc123} demonstrate the process of sending $\tau$ through the classically-$123$-avoiding stack $s_{123}=\sc_{123}$, the consecutively-$123$-avoiding stack $\sc_{\un{123}}$, and the vincular-avoiding-stacks $\sc_{\un{12}3}$ and $\sc_{1\un{23}}$, respectively. We find that $\sc_{123}(\tau)=463215$, $\sc_{\un{123}}(\tau)=263415$, $\sc_{\un{12}3}(\tau)=426315$, and $\sc_{1\un{23}}(\tau)=632415$. In particular, the four maps produce different images after acting on $\tau$. 

\begin{figure}[h]
    \centering
    \includegraphics[width=\textwidth]{./pictures/exs123}
    \caption{The stack-sorting map $s_{123}=\sc_{123}$ acting on $\tau=514362$}
    \label{exs123}
\end{figure}

\begin{figure}[h]
    \centering
    \includegraphics[width=\textwidth]{./pictures/exsc123}
    \caption{The stack-sorting map $\sc_{\un{123}}$ acting on $\tau=514362$}
    \label{exsc123}
\end{figure}

\begin{figure}[h]
    \centering
    \includegraphics[width=\textwidth]{./pictures/exscs123}
    \caption{The stack-sorting map $\sc_{\un{12}3}$ acting on $\tau=514362$}
    \label{exscs123}
\end{figure}

\begin{figure}[h]
    \centering
    \includegraphics[width=\textwidth]{./pictures/exssc123}
    \caption{The stack-sorting map $\sc_{1\un{23}}$ acting on $\tau=514362$}
    \label{exssc123}
\end{figure}

\medskip

It is well-known \cite{Knuth} that the permutations sorted onto the identity by $s$ are those that avoid the pattern $231$. Since \citet{West}'s introduction of the deterministic stack-sorting map, researchers studied the permutations sorted onto the identity by several applications of $s$ \cite{BonaStackSortingSurvey, DefantPreimagesStackSorting, DefantThreeStackSortable, Knuth, West, ZeilbergerTwoStackSortable}. The map $\sc_\s$ has also been studied extensively from a sorting point of view. For classical and consecutive patterns $\s$, the articles \cite{TwoStacks, CCF, 132machine, Defant} focused on the permutations $\sort(\sc_\s)$ sorted to the identity by $s\circ\sc_\s$. \citet{CCF}, which first introduced and studied the map $\sc_\s$ for classical patterns $\s$, was able to describe and enumerate the sorting class $\sort(\sc_\s)$ for $\s=123$ and $\s=321$. Since then, variations and generalizations of \citet{CCF}'s work have been investigated \cite{TwoStacks, 132machine, Defant}. The sorting class of $\sc_{132}$ was found by \cite{132machine}, and the sorting classes of $\sc_{\un{132}}$, $\sc_{\un{123}}$, and $\sc_{\un{321}}$ were found by \citet{Defant}. 

\medskip

Another recent topic of study is the dynamical properties of $\sc_\s$, first examined by \citet{KatalinBerlow}. Various authors worked on the maximum number of preimages under $s$ and $\sc_\s$ \cite{KatalinBerlow, DefantPreimagesStackSorting, DefantFertility, Defant}. Most relevantly to our work, \citet{KatalinBerlow} showed that, for patterns $\s$ of length $3$, the maximum number of preimages a permutation can have under $\sc_\s$ is given by $C_{n-1}$. Furthermore, the periodic points of the map $\sc_\s$ have been studied \cite{KatalinBerlow, ChoiChoi, Defant}, especially for classical and consecutive patterns $\s$ of length $3$. In particular, \citet{Defant} showed that the periodic points of $\sc_{132}$ are the permutations that avoid the patterns $132$ and $231$, and the periodic points of $\sc_{312}$ are the permutations that avoid the patterns $312$ and $213$. 

\medskip

The remainder of the paper is organized as follows. In \Cref{sortingclass}, we describe and enumerate $\sort(\sc_\s)$ for various vincular patterns $\s$. Additionally, we decide whether $\sort(\sc_\s)$ is a permutation class. Our results are summarized in \Cref{sortingclassresults}. Henceforth, $C_n$, $M_n$, and $S_n$ denote the $n$-th Catalan, Motzkin, and Schr\"oder number, respectively. 

\medskip

\begin{footnotesize}
\begin{table}[h]
\centering
\begin{tabular}{c | c c c c c c c c c c} 
    $\s$ & $\sort_n(\sc_\s)$ & $|\sort_n(\sc_\s)|$ & OEIS & Permutation class \\ 
    \hline \hline
    1\un{23} & $\av_n(3214,4213,132)$ & $2^n-n$ & \href{https://oeis.org/A000325}{A000325} & Yes \\ 
    \un{32}1 & $\av_n(\un{123},132)$ & $M_n$ & \href{https://oeis.org/A001006}{A001006} & No \\
    3\un{21} & $\av_n(123,132)$ & $2^{n-1}$ & \href{https://oeis.org/A000079}{A000079} & Yes \\
    \un{13}2 & $\av_n(2314,\mu_{132})$ & $\sum_{i=0}^{n-1}\binom{n-1}iC_i$ & \href{https://oeis.org/A007317}{A007317} & No \\ 
    \un{23}1 & $\av_n(1324,\mu_{2413})$ & $S_{n-1}$ & \href{https://oeis.org/A006318}{A006318} & No \\ 
    1\un{32} & \Cref{ssc132characterization} & $C_n$ & \href{https://oeis.org/A000108}{A000108} & No \\ 
    \un{12}3 & \Cref{scs123characterization} & $C_{n-1}+\sum_{i=0}^{n-2}2^{n-2-i}C_i$ & \href{https://oeis.org/A126221}{A126221} & No
\end{tabular}
\caption{The sorting classes $\sort_n(\sc_\s)$}
\label{sortingclassresults}
\end{table}
\end{footnotesize}

\medskip

We explore the dynamics of $\sc_\s$ in \Cref{dynamike}. In \Cref{preimages}, we show that the $C_{n-1}$ bound proposed by \citet{KatalinBerlow} is tight for $\s=\un{12}3$, $\un{32}1$, $\un{23}1$, $\un{21}3$, $2\un{31}$, and $2\un{13}$. In addition, we conjecture and provide progress towards showing that \[\max_{\pi\in \sy_n}|\sc_{1\un{23}}^{-1}(\pi)|=\max_{\pi\in \sy_n}|\sc_{3\un{21}}^{-1}(\pi)|=2^{n-2}.\]
In \Cref{periodicpoints}, we briefly consider the periodic points of $\sc_\s$. We show in particular that the periodic points of $\sc_\s$ for $\s\in\{\un{13}2, \un{23}1, 1\un{32}, 2\un{31}\}$ are the permutations that avoid the patterns $132$ and $231$, and for $\s\in\{\un{31}2, \un{21}3, 3\un{12},2\un{13}\}$ are the permutations that avoid the patterns $312$ and $213$. 

\medskip

\section{Preliminaries} \label{prelim}
\subsection{Definitions and Notation} \label{definitions}

Consider $\pi=\pi_1\pi_2\cdots\pi_n$. We define the \emph{reverse} of $\pi$ as $\rev(\pi)=\pi_n\pi_{n-1}\cdots\pi_1$ and the \emph{complement} of $\pi$ as $\comp(\pi)=(n+1-\pi_1)(n+1-\pi_2)\cdots(n+1-\pi_n)$. For instance, $\rev(1342)=2431$ and $\comp(1342)=4213$. These definitions extend easily to patterns $\s$. An element $\tau_i$ is a \emph{left-to-right minimum} of $\tau$ if, for all $j<i$, we have $\tau_i<\tau_j$. We often split a permutation $\tau$ into the concatenation of two or more subsequences. For example, we may let $\tau=\tau^a\tau^b$, where $\tau^a=\tau_1\cdots\tau_i$ and $\tau^b=\tau_{i+1}\cdots\tau_n$. We consider $a$ and $b$ as permutations in their own right, and they can be processed by stack-sorting maps. 

\medskip

Two sequences $\tau$ and $\pi$ of size $n$ are \emph{order-isomorphic} if, for each pair of indices $i$ and $j$, we have $\tau_i<\tau_j$ if and only if $\pi_i<\pi_j$. We write $\tau\simeq\pi$ if $\tau$ and $\pi$ are order-isomorphic. For a permutation $\pi$ and pattern $\s$, we say $\pi$ \emph{avoids} $\s$ if there exists no subsequence of $\pi$ order-isomorphic to $\s$. Else, we say $\pi$ \emph{contains} $\s$. As a generalization, we allow $\s$ to be a vincular pattern, where terms in the subsequence of $\pi$ corresponding to underlined terms in $\s$ must be adjacent. Henceforth, we work only with vincular pattern avoidance. Let $\S_n$ be the set of the vincular patterns of length $n$, and let $\S=\bigcup_{n\ge 1}\S_n$.

\medskip

If $\Sigma$ is a list of patterns, let $\av_n(\Sigma)$ be the permutations in $\sy_n$ that avoid all the patterns in $\Sigma$. Let $\av(\Sigma)=\bigcup_{n\ge 1}\av_n(\Sigma)$. It is well-known \cite{WestCatalan} that $|\av_n(\s)|=C_n$ for each $\s\in S_3$. 

\medskip

Throughout the remainder of the article, let $\sc_\s$ denote the stack-sorting map which uses a stack that avoids $\s\in\S$. Note that $s=\sc_{21}$. Naturally, the stack-sorting process of $\sc_\s$ keeps track of three sequences: the input sequence, the stack, and the output sequence. For brevity, we describe an element as being popped from or added to the stack of $\sc_\s$ without explicitly mentioning the stack. A \emph{premature} element $\tau_i$ is one that is popped before the input sequence is empty. For example, the elements $4$, $6$, and $3$ are premature in \Cref{exs123}. Note that $\tau_1$ is not premature when $\tau$ is sorted by $\sc_\s$ for $\s\in \S_3$. 

\medskip

\subsection{Preliminary Lemmas}

\Cref{complemma} follows from the definition of the map $\sc_\s$. 

\medskip

\begin{lemma}[{\citet[Lemma 1.2]{Defant}}] \label{complemma}
    For a pattern $\s\in \S$ of length at least $2$, we have 
    \[\sc_{\comp(\s)}=\comp\circ\sc_\s\circ\comp.\]
\end{lemma}

\medskip

Next, we show that avoiding some vincular patterns is equivalent to avoiding others. 

\medskip

\begin{lemma} \label{free2avoids}
    For pairs $(\s_1,\s_2)\in\{(\un{13}2,132),(\un{31}2,312),(2\un{13},213),(2\un{31},231)\}$, we have $\av(\s_1)=\av(\s_2)$. 
\end{lemma}

\begin{proof}
    We show that containing $\un{13}2$ is equivalent to containing $132$, where the other cases follow similarly. If $\tau\in \sy_n$ contains $\un{13}2$, then it clearly also contains $132$. For the other direction, suppose $\tau$ contains $\tau_i\tau_j\tau_k\simeq 132$. But there must exist some index $i\le l<j$ satisfying $\tau_l<\tau_k<\tau_{l+1}$, so that $\tau$ contains $\tau_l\tau_{l+1}\tau_k\simeq \un{13}2$. 
\end{proof}

\medskip

\Cref{free2avoids} implies the following.

\medskip

\begin{corollary} \label{free2stacks}
    The maps $\sc_{\un{13}2}$, $\sc_{\un{31}2}$, $\sc_{2\un{31}}$, and $\sc_{2\un{13}}$ are equivalent to $\sc_{132}$, $\sc_{312}$, $\sc_{231}$, and $\sc_{213}$, respectively. 
\end{corollary}

\medskip

\section{The sorting classes $\sort(\sc_\s)$} \label{sortingclass}
In this section, we study the \emph{sorting class} $\sort_n(\sc_\s)$ of $\sc_\s$, the set of permutations in $\sy_n$ sorted onto the identity by $s\circ\sc_\s$. Let $\sort(\sc_\s)=\bigcup_{n\ge 1}\sort_n(\sc_\s)$. Additionally, a \emph{permutation class} is a (possibly infinite) set $\Pi$ of permutations such that every $\tau$ contained in $\pi\in \Pi$ is also in $\Pi$. We begin with the following well-known lemma. 

\medskip

\begin{lemma}[\citet{Knuth}]
    The set of permutations sorted by $s$ to the identity is $\av(231)$. 
\end{lemma}

\medskip

It thus suffices to find the permutations sent by $\sc_\s$ to $\av(231)$. We now characterize and enumerate $\sort_n(\sc_\s)$ for several $\s\in\S_3$. Furthermore, we decide whether $\sort(\sc_\s)$ is a permutation class. 

\medskip

\subsection{The Sorting Classes of $\sc_{1\un{23}}$, $\sc_{\un{32}1}$, and $\sc_{3\un{21}}$} \label{avoidance}

In this section, we study sorting classes that can be described by pattern avoidance. The characterizations and sizes of $\sort_n(\sc_{1\un{23}})$, $\sort_n(\sc_{\un{32}1})$, and $\sort_n(\sc_{3\un{21}})$ are listed in \Cref{sortingclassresults}. We begin by describing one aspect of $\sc_{1\un{23}}$ in \Cref{1un23lemma}. 

\medskip

\begin{lemma}\label{1un23lemma}
    Consider $\tau\in \sy_n$ that is being processed by $\sc_{1\un{23}}$, and suppose $\tau_i$ is an element currently in the stack. Subsequently, after each time an element enters the stack, the stack will contain some element of value at most $\tau_i$. 
\end{lemma}
\begin{proof}
    Let $\tau_p$ be the next element to be added to the stack. If no elements are popped immediately before $\tau_p$ is added, then $\tau_i$ remains in the stack. Else, suppose $\tau_q$ is the element popped immediately before $\tau_p$ is added to the stack so that $\tau_p<\tau_q$. Then, $\tau_q\le \tau_i$, for otherwise the stack would already contain $1\un{23}$. Thus, $\tau_p<\tau_i$, and the conclusion follows. 
\end{proof}

\medskip

With \Cref{1un23lemma}, we are ready to directly describe $\sort_n(\sc_{1\un{23}})$. 

\medskip

\begin{theorem}\label{ssc123characterization}
    The sorting class of $\sc_{1\un{23}}$ is characterized by \[\sort_n(\sc_{1\un{23}})=\av_n(132,3214,4213).\]
\end{theorem}
\begin{proof}
    Let $\pi=\sc_{1\un{23}}(\tau)$. We first show that if $\tau$ does not avoid $132$, $3214$, and $4213$, then $\pi$ contains $231$. If $\tau$ contains $\tau_i\tau_j\tau_k\simeq 132$, \Cref{1un23lemma} implies that when $\tau_j$ enters the stack, there is some term $\tau_i'\le \tau_i$ below it. We split into three possibilities: 
    \begin{itemize}
        \item If $\tau_j$ has not popped when $\tau_k$ enters (and neither has $\tau_i'$), then the stack, and thus $\pi$ will contain $\tau_k\tau_j\tau_i'\simeq231$. 
        \item Suppose $\tau_j$, but not $\tau_i'$, was popped before $\tau_k$ enters. Suppose further that $\tau_j$ but not $\tau_r$ was popped when $\tau_p$ was being added to the stack, because $\tau_p\tau_q\tau_r\simeq1\un{23}$. If $\tau_r\ne\tau_i'$ and $\tau_i'<\tau_q$, then $\tau_q\tau_r\tau_i'\simeq 231$ within the stack and thus $\pi$. As such, it suffices to consider only the case $\tau_q<\tau_l$. Meanwhile, if $\tau_r=\tau_i'$, then $\tau_q<\tau_r=\tau_i'<\tau_k$. We proceed with the assumption $\tau_q<\tau_k$. By \Cref{1un23lemma}, when $\tau_k$ enters the stack, there is some term $\tau_p'\le \tau_p$ below it. Therefore, $\pi$ contains $\tau_q\tau_k\tau_p'\simeq231$. 
        \item Lastly, suppose that $\tau_i'$ was popped before $\tau_k$ enters the stack. As above, consider $\tau_p$, $\tau_q$, and $\tau_r$. If $\tau_q\ne \tau_i'$, we must have $\tau_p<\tau_i'$, as otherwise the stack contains $\tau_i'\tau_q\tau_r\simeq1\un{23}$. Otherwise, if $\tau_q=\tau_i'$, then $\tau_p<\tau_q=\tau_i'$. In either case, we have $\tau_p<\tau_i'$. Again by \Cref{1un23lemma}, when $\tau_k$ is at the top of the stack, there is some term $\tau_p'\le \tau_p$ below it. Therefore, $\pi$ contains $\tau_i'\tau_k\tau_p'\simeq 231$. 
    \end{itemize}

    \medskip
    
    As per the above cases, $\tau$ will not be in the sorting class of $\sc_{1\un{23}}$ if it contains $132$. Henceforth, assume that $\tau\in\av_n(132)$. If $\tau$ contains $\tau_i\tau_j\tau_k\tau_l\simeq3214$ and neither of $\tau_i$ and $\tau_j$ have popped when $\tau_k$ enters the stack, then the stack will contain $\tau_k\tau_j\tau_i\simeq 123$. In order for the stack to avoid the pattern $1\un{23}$, there must be some element $\tau_s<\tau_k$ between $\tau_i$ and $\tau_j$ in the stack. But then $\tau$ contains $\tau_s\tau_j\tau_k\simeq 132$, a contradiction. Thus, at least one of $\tau_i$ or $\tau_j$ will be popped before $\tau_k$ enters. By \Cref{1un23lemma}, there is some $\tau_k'\le \tau_k$ in the stack when $\tau_l$ is added. Then, $\pi$ will contain either $\tau_i\tau_l\tau_k'\simeq231$ or $\tau_j\tau_l\tau_k'\simeq231$. Lastly, if $\tau$ contains $\tau_i\tau_j\tau_k\tau_l\simeq4213$, we proceed similarly as before in most cases. The only new scenario we must consider is if $\tau_i$, but not $\tau_j$, is popped when $\tau_k$ enters, as $\tau_i\tau_l\tau_k'\not\simeq231$. However, note that $\tau_l<\tau_1$, for otherwise $\tau_1\tau_i\tau_l\simeq132$. So, we set $i=1$ without changing the relative ordering of $\tau_i\tau_j\tau_k\tau_l$. But $\tau_1$ can never be popped, so we may simply proceed as before. In conclusion, if $\tau$ is not in $\av_n(132,3214,4213)$, it is not in $\sort_n(\sc_{1\un{23}})$. 

    \medskip

    Next, we show that if $\pi$ contains $\pi_i\pi_j\pi_k\simeq 231$, then $\tau$ contains at least one of $132$, $3214$, or $4213$. We split into the following cases based on the order in which $\pi_i$, $\pi_j$, and $\pi_k$ appear in $\tau$: 
\begin{itemize}
    \item Suppose $\pi_i$ occurs before $\pi_j$ in $\tau$, so that $\pi_i$ must pop before $\pi_j$ enters the stack. Say $\pi_i$ was popped when $\tau_p$ was being added, as $\tau_p\tau_q\tau_r\simeq 1\un{23}$. Note that $\tau_q\le \pi_i$, for otherwise the stack would contain $\pi_i\tau_q\tau_r\simeq 1\un{23}$. Then, $\tau$ contains $\tau_r\tau_q\tau_p\pi_j$, which forms either the $3214$ or $4213$ pattern. 
    \item If $\tau$ contains $\pi_k\pi_j\pi_i\simeq 132$, we are done. 
    \item Next, suppose $\tau$ contains $\pi_j\pi_i\pi_k$, so that $\pi_i$ pops before $\pi_k$ enters the stack. As before, say $\pi_i$ was popped when $\tau_p$ was being added to the stack, as $\tau_p\tau_q\tau_r\simeq 1\un{23}$. If $\tau_q<\pi_i$, then $\tau$ contains $\tau_q\pi_j\pi_i\simeq 132$. Else, $\pi_i<\tau_q$, so that $\pi_i\tau_q\tau_r\simeq 1\un{23}$, but then $\pi_j$ would have popped before $\pi_i$ entered, a contradiction. 
    \item Lastly, if $\tau$ contains $\pi_j\pi_k\pi_i$, the element $\pi_j$ must pop before $\pi_k$ enters but after $\pi_i$ enters, which is impossible. 
\end{itemize}
\end{proof}

\medskip

\begin{proposition}
    It holds that $\sort(\sc_{1\un{23}})$ is a permutation class as every permutation in $\sort(\sc_{1\un{23}})$ necessarily avoids every permutation outside of $\sort(\sc_{1\un{23}})$. 
\end{proposition}

\medskip

The following two lemmas allow us to enumerate $\av_n(132,3214,4213)$ in \Cref{av13232144213enumeration}. 

\medskip

\begin{lemma}[{\citet[Proposition 11]{SimionSchmidt}}] \label{132321enumeration} 
    It holds that $|\av_n(132,321)|=\binom n2+1$.
\end{lemma}

\medskip

\begin{lemma}[{\citet[Proposition 8]{SimionSchmidt}}] \label{132213enumeration} 
    It holds that $|\av_n(132,213)|=2^{n-1}$ for $n\ge 1$. Additionally, $|\av_0(132,213)|=1$. 
\end{lemma}

\medskip

\begin{theorem} \label{av13232144213enumeration}
    It holds that $|\av_n(132,3214,4213)|=2^n-n$. 
\end{theorem}
\begin{proof}
    Let $\tau^a$ and $\tau^b$ be subsequences of $\tau$ so that $\tau=\tau^an\tau^b$. We claim that the following three conditions are equivalent to $\tau\in\av_n(132,3214,4213)$.  
    \begin{enumerate}
        \item each element of $\tau^a$ is larger than each element of $\tau^b$, 
        \item $\tau^a\in\av(132,321)$, and
        \item $\tau^b\in\av(132,213)$. 
    \end{enumerate}

    \medskip
    
    Condition (1) follows by letting $n$ correspond to the $3$ in the pattern $132$. Conditions (2) and (3) follow by letting $n$ correspond to the $4$ in the patterns $3214$ and $4213$. To check that the conditions are sufficient, suppose they all hold. Let $i$ be the index such that $\tau_i=n$. 
    
    \begin{itemize}
        \item We show that each subsequence $\tau_j\tau_k\tau_l$ of $\tau$ avoids $132$. Clearly, we have $j\ne i$ and $l\ne i$. Since $\tau^a$ and $\tau^b$ individually avoid $132$, we only need to consider when $j<i<l$. But it suffices to let $k=i$, which is taken care of by condition (1). 
        \item We show that each subsequence $\tau_j\tau_k\tau_l\tau_m$ of $\tau$ avoids $3214$. Clearly, we must have $j,k,l\ne i$. If $m\le i$, it suffices to only consider $m=i$, which is satisfied as $\tau^a$ avoids $321$. Else, if $i<m$, the inequality $\tau_j<\tau_m$ along with the condition (1) implies that $i<j$. But $\tau^b$ avoids $213$, and thus also $3214$. 
        \item We show that each subsequence $\tau_j\tau_k\tau_l\tau_m$ of $\tau$ avoids $4213$. Clearly, we must have $k,l,m\ne i$. If $i\le j$ or if $j\le i<k$, only need to consider $j=i$. Since $\tau^b$ avoids $213$, the case $j=i$ is covered. Note that the cases $l<i<k$ and $i<k<l$ are impossible by condition (1) and the inequality $\tau_l<\tau_k<\tau_m$. Lastly, say $m<i$. But $\tau^a$ avoids $321$, and thus also $4213$. 
    \end{itemize}

    \medskip

    Due to the condition (1), the index $i$ determines the values of the elements of $\tau^a$ and $\tau^b$. Those elements are then independently ordered by the conditions (2) and (3). We thus sum over $i$, using \Cref{132321enumeration,132213enumeration}:
    \begin{align*}
        |\av_n(132,3214,4213)|&=\sum_{i=1}^n |\av_{i-1}(132,321)|\cdot |\av_{n-i}(132,213)|\\
        &=\left(\binom{n-1}2+1\right)+\sum_{k=0}^{n-2}\left(\binom k2+1\right)2^{n-k-2}\\ 
        &=\left(\binom{n-1}2+1\right)+\sum_{k=0}^{n-2}\binom k2 2^{n-2-k}+\sum_{k=0}^{n-2}2^{n-k-2}\\
        &=\left(\frac12(n-1)(n-2)+1\right)+\left(2^{n-1}-1-\frac12n(n-1)\right)+(2^{n-1}-1)\\
        &=2^n-n,
    \end{align*}
    as desired. 
\end{proof}

\medskip

\Cref{ssc123enumeration} then follows from \Cref{ssc123characterization} and \Cref{av13232144213enumeration}.

\medskip

\begin{corollary}\label{ssc123enumeration}
    The sorting class of $\sc_{1\un{23}}$ is enumerated by \[|\sort_n(\sc_{1\un{23}})|=2^n-n.\]
\end{corollary}

\medskip

\Cref{ssc123enumeration} completes our analysis of $\sort_n(\sc_{1\un{23}})$. Next, we find $\sort_n(\sc_{\un{32}1})$ and $\sort_n(\sc_{3\un{21}})$. Recall the following results, which we generalize in \Cref{321general}: 

\medskip

\begin{theorem}[{\citet[Corollary 3.3]{CCF}}] \label{321CCF} 
    For $k\ge 3$, we have \[\sort_n(\sc_{k(k-1)\cdots1})=\av_n(12\cdots k,132).\]
\end{theorem}

\medskip

\begin{theorem}[{\citet[Theorem 6.1]{Defant}}] \label{321Defant}
    For $k\ge 3$, we have \[\sort_n(\sc_{\un{k(k-1)\cdots1}})=\av_n(\un{12\cdots k},132).\]
\end{theorem}

\medskip

The following is a natural extension of \Cref{321CCF,321Defant}. 

\medskip

\begin{theorem} \label{321general}
    For $k\ge 3$, let $\s\in \S_k$ be a pattern with the values $k(k-1)\cdots 1$. Then, \[\sort_n(\sc_\s)=\av_n(\rev(\s),132).\]
\end{theorem}
\begin{proof}
    Let $\pi=\sc_\s(\tau)$. Suppose that $\tau$ contains $\rev(\s)$, and let the leftmost occurrence of $\rev(\s)$ in $\tau$ be $\tau_{l_1}\tau_{l_2}\cdots\tau_{l_k}$. Evidently, $\pi$ contains $\tau_{l_{k-1}}\tau_{l_k}\tau_{l_{k-2}}\simeq231$. If $\tau$ avoids $\rev(\s)$, then $\sc_\s(\tau)=\rev(\tau)$. As such, $\tau$ must also avoid $\rev(231)=132$. 

    \medskip

    To check the converse, suppose $\pi$ contains $231$. Then, either $\pi\ne \rev(\tau)$ or $\pi=\rev(\tau)$. In the first scenario, $\tau$ must contain $\rev(\s)$. In the second, $\tau$ must contain $132$.     
\end{proof}

\medskip

It follows from \Cref{321general} the characterizations of $\sort_n(\sc_{\un{32}1})$ and $\sort_n(\sc_{3\un{21}})$. \Cref{scs321simplification} then simplifies our description of $\sort_n(\sc_{\un{32}1})$. 

\medskip

\begin{lemma} \label{scs321simplification}
    It holds that $\av_n(1\un{23},132)=\av_n(\un{123},132)$.
\end{lemma}
\begin{proof}
        Clearly, each permutation in $\av_n(1\un{23},132)$ is also in $\av_n(\un{123},132)$. Now, suppose for contradiction that $\tau$ avoids $\un{123}$ and $132$ but not $1\un{23}$. Say $\tau$ contains $\tau_{i}\tau_{j-1}\tau_{j}\simeq1\un{23}$, with $i<j-2$. Since $\tau$ avoids $\un{123}$, $\tau_{j-2}$ is required to be greater than $\tau_{j-1}$. But then $\tau_{i}\tau_{j-2}\tau_{j-1}$ is an instance of $132$. The conclusion follows. 
\end{proof}

\medskip

\Cref{321general} and \Cref{scs321simplification} give the following characterization of $\sort_n(\sc_{\un{32}1})$. 

\medskip

\begin{corollary} \label{scs321characterization}
    The sorting class of $\sc_{\un{32}1}$ is characterized by  \[\sort_n(\sc_{\un{32}1})=\av_n(\un{123},132).\] 
\end{corollary}

\medskip

\begin{proposition}
    It does not hold that $\sort(\sc_{\un{32}1})$ is a permutation class since $2314$ contains $123$ but $\sc_{\un{32}1}(2314)=4132$ avoids $231$ while $\sc_{\un{32}1}(123)=231$ does not avoid $231$. 
\end{proposition}

Note that from \Cref{321Defant}, we have $\sort_n(\sc_{\un{321}})=\av_n(\un{123},132)=\sort_n(\sc_{\un{32}1})$. Additionally, for $\tau\in \av_n(\un{123},132)$, we note that $\sc_{\un{321}}(\tau)=\rev(\tau)=\sc_{\un{32}1}(\tau)$. However, the two maps are distinct. For example, $\sc_{\un{321}}(1324)=4231\ne 2341=\sc_{\un{32}1}(1324)$. Next, \Cref{un123132enumeration} allows us to enumerate $\sort_n(\sc_{\un{32}1})$.

\medskip

\begin{lemma}[{\citet[Proposition 3.2]{ModularCatalan}}] \label{un123132enumeration} 
    It holds that $|\av_n(\un{12\cdots k},132)|=M_{k-1,n}$, where $M_{m,n}$ denotes the generalized Motzkin number, for $k\ge 3$. 
\end{lemma}

\medskip

 Recalling that $M_{2,n}=M_n$ \cite{ModularCatalan}, \Cref{scs321enumeration} follows from \Cref{scs321characterization} and \Cref{un123132enumeration}. 

\medskip

 \begin{corollary} \label{scs321enumeration}
     The sorting class of $\sc_{\un{32}1}$ is enumerated by  \[|\sort_n(\sc_{\un{32}1})|=M_n.\]
 \end{corollary}

\medskip

Having completed our analysis of $\sort_n(\sc_{\un{32}1})$, we move on to characterize $\sort_n(\sc_{3\un{21}})$. \Cref{ssc321simplification} simplifies our description of $\sort_n(\sc_{3\un{21}})$. 

\medskip

\begin{lemma} \label{ssc321simplification}
    It holds that $\av_n(\un{12}3,132)=\av_n(123,132)$.
\end{lemma}
\begin{proof}
    Clearly, each permutation in $\av_n(123,132)$ is also in $\av_n(\un{12}3,132)$. Next, suppose for contradiction that $\tau$ avoids $\un{12}3$ and $132$ but not $123$. Say $\tau$ contains $\tau_i\tau_j\tau_k\simeq 123$, where $i<j-1$. If $\tau_{j-1}<\tau_j$, then $\tau_{j-1}\tau_j\tau_k\simeq \un{12}3$. Else, if $\tau_j<\tau_{j-1}$, then $\tau_i\tau_{j-1}\tau_j\simeq132$. The conclusion follows.
\end{proof}

\medskip

\Cref{321general} and \Cref{ssc321simplification} then imply \Cref{ssc321characterization}. 

\medskip

\begin{corollary} \label{ssc321characterization}
    The sorting class of $\sc_{3\un{21}}$ is characterized by  \[\sort_n(\sc_{3\un{21}})=\av_n(123,132).\]
\end{corollary} 

\medskip

\begin{proposition}
    It holds that $\sort(\sc_{3\un{21}})$ is a permutation class as every permutation in $\sort(\sc_{3\un{21}})$ necessarily avoids every permutation outside of $\sort(\sc_{3\un{21}})$.
\end{proposition}

\medskip

From \Cref{321CCF}, we have $\sort_n(\sc_{321})=\av_n(123,132)=\sort_n(\sc_{3\un{21}})$. Furthermore, for $\tau\in \av_n(123,132)$, we note that $\sc_{321}(\tau)=\rev(\tau)=\sc_{3\un{21}}(\tau)$. However, the two maps are distinct. For example, $\sc_{321}(1324)=4231\ne 2341=\sc_{3\un{21}}(1324)$. Next, \Cref{123132enumeration} allows us to enumerate $\sort_n(\sc_{3\un{21}})$. 

\medskip

\begin{lemma}[{\citet[Lemma 5]{SimionSchmidt}}] \label{123132enumeration} 
    It holds that $|\av_n(123,132)|=2^{n-1}$ for $n\ge 1$. Additionally, $|\av_0(123,132)|=1$. 
\end{lemma}

\medskip

\Cref{ssc321enumeration} follows from \Cref{ssc321characterization} and \Cref{123132enumeration}, completing our analysis of $\sort_n(\sc_{3\un{21}})$. 

\medskip

\begin{corollary} \label{ssc321enumeration}
    The sorting class of $\sc_{3\un{21}}$ is enumerated by  \[|\sort_n(\sc_{3\un{21}})|=2^{n-1}.\]
\end{corollary}

\medskip

\subsection{The Sorting Classes of $\sc_{\un{13}2}$ and $\sc_{\un{23}1}$} \label{mesh}

\begin{figure}[h]
    \centering
    \includegraphics[width=8cm]{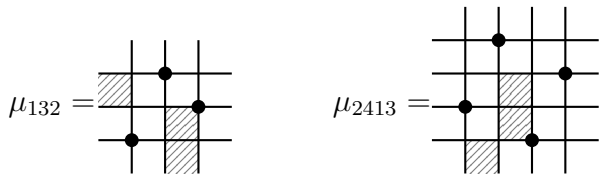}
    \caption{The mesh patterns $\mu_{132}$ and $\mu_{2413}$.}
    \label{meshdiagram}
\end{figure}

In this section, we introduce \emph{mesh patterns} in order to describe $\sort_n(\sc_{\un{13}2})$ and $\sort_n(\sc_{\un{23}1})$. A mesh pattern $\mu$ of length $k$ consists of a pattern $\s\in S_k$ and a set $A\in [k]^2$ of forbidden squares in the plot of $\s$, identified by their lower left corners. A permutation $\pi$ contains $\mu=(\s,A)$ if $\pi$ has an occurrence of $\s$ such that no elements of $\pi$ lie within the forbidden squares in $A$. Henceforth, let $\mu_{132}=(132,\{(0,2),(2,0),(2,1)\})$ and let $\mu_{2413}=(2413,\{(1,0),(2,1),(2,2)\})$. The plots of $\mu_{132}$ and $\mu_{2413}$ are shown in \Cref{meshdiagram}. 

We will show that $\sort_n(\sc_{\un{13}2})$ and $\sort_n(\sc_{\un{23}1})$ can be described by mesh pattern avoidance. The characterizations of the sorting classes are listed in \Cref{sortingclassresults}, along with their sizes. 

\medskip

We begin with $\sort_n(\sc_{\un{13}2})$. Recall that, by \Cref{free2stacks}, we have $\sc_{\un{13}2}=\sc_{132}$. \citet{132machine} identified the class $\sort_n(\sc_{132})$ with mesh pattern avoidance, and enumerated the class using a bijection between $\sort_n(\sc_{132})$ and the restricted growth functions of length $n$ that avoid $12231$. 

\medskip

\begin{theorem}[{\citet[Theorem 3.4]{132machine}}] \label{132characterization} 
    The sorting class of $\sc_{132}$ is characterized by \[\sort_n(\sc_{132})=\av_n(2314,\mu_{132}).\] 
\end{theorem}

\medskip

\begin{theorem}[{\citet[Corollary 4.10]{132machine}}] \label{132enumeration} 
    The sorting class of $\sc_{132}$ is enumerated by \[|\sort_n(\sc_{132})|=\sum_{i=0}^{n-1}\binom{n-1}iC_i.\] 
\end{theorem}

\medskip

As immediate corollaries of \Cref{132characterization,132enumeration}, we describe $\sort(\sc_{\un{13}2})$. 

\medskip

\begin{corollary} \label{scs132characterization}
    The sorting class of $\sc_{\un{13}2}$ is characterized by \[\sort_n(\sc_{\un{13}2})=\av_n(2314,\mu_{132}).\] 
\end{corollary}

\medskip

\begin{proposition}
    It does not hold that $\sort(\sc_{\un{13}2})$ is a permutation class since $2413$ contains $132$ but $\sc_{\un{13}2}(2413)=4312$ avoids $231$ while $\sc_{\un{13}2}(132)=231$ does not avoid $231$. 
\end{proposition}

\medskip

\begin{corollary} \label{scs132enumeration}
    The sorting class of $\sc_{\un{13}2}$ is enumerated by \[|\sort_n(\sc_{\un{13}2})|=\sum_{i=0}^{n-1}\binom{n-1}iC_i.\] 
\end{corollary}

\medskip

Next, we show that $\sort_n(\sc_{\un{23}1})$ can also be described with mesh pattern avoidance. We first show \Cref{un231lemma} about one property of $\sc_{\un{23}1}$. 

\medskip

\begin{lemma}    \label{un231lemma}
    Consider $\tau\in \sy_n$ that is being processed by $\sc_{\un{23}1}$, and suppose $\tau_i$ is an element currently in the stack. Subsequently, after each time an element enters the stack, the stack will contain some element of value at most $\tau_i$. 
\end{lemma}
\begin{proof}
    When $\tau_i$ is popped from the stack, the next element to be added to the stack must be less than $\tau_i$. 
\end{proof}

\medskip

With \Cref{un231lemma}, we describe $\sort_n(\sc_{\un{23}1})$ in \Cref{scs231characterization}. 

\medskip

\begin{theorem} \label{scs231characterization}
    The sorting class of $\sc_{\un{23}1}$ is characterized by \[\sort_n(\sc_{\un{23}1})=\av_n(1324,\mu_{2413}).\] 
\end{theorem}
\begin{proof}
    Let $\pi=\sc_{\un{23}1}(\tau)$. We first show that if $\tau$ does not avoid $1324$ and $\mu_{2413}$, then $\pi$ contains $231$. Suppose $\tau$ contains $\tau_i\tau_j\tau_k\tau_l\simeq 1324$. By \Cref{un231lemma}, there is some term $\tau_i'\le \tau_i$ in the stack when $\tau_j$ enters. If $\tau_j$ is not popped when $\tau_k$ enters, the stack, and thus $\pi$, would contain $\tau_k\tau_j\tau_i'\simeq 231$. Thus, $\tau_j$ is popped before $\tau_k$ enters. \Cref{un231lemma} implies that there is some $\tau_k'\le \tau_k$ in the stack when $\tau_l$ enters. Then, $\pi$ contains $\tau_j\tau_l\tau_k'\simeq 231$. Now assuming that $\tau$ avoids $1324$, say $\tau$ contains $\tau_i\tau_j\tau_k\tau_l\simeq \mu_{2413}$. In what follows, we consider the elements of $\tau$ relative to the plot of $\tau_i\tau_j\tau_k\tau_l$. Note that the $(0,0)$ square is forbidden, as any element in that square would form the $1324$ pattern in $\tau$ along with $\tau_i\tau_k\tau_l$. By \Cref{un231lemma}, there is some element $\tau_i'\le \tau_i$ in the stack when $\tau_j$ enters. Since the $(1,0)$ square is banned, we have that $\tau_k\le \tau_i'$. We now claim that there will be some element $\tau_j'\ge \tau_l$ in the stack when $\tau_l$ enters. Suppose $\tau_j$ was popped before $\tau_k$ enters the stack while $\tau_p$ was being added because $\tau_p\tau_j\tau_r\simeq \un{23}1$. Then, $\tau_k<\tau_r$, which implies that $\tau_k<\tau_p$. As $(2,1)$ and $(2,2)$ are forbidden, we have that $\tau_l<\tau_r$, which implies that there will be some $\tau_j'\ge \tau_l$ in the stack when $\tau_k$ enters. But after $\tau_k$ enters, $\tau_j'$ cannot be popped, so that the stack, and thus $\pi$, contains $\tau_l\tau_j'\tau_i'\simeq 231$. 

    \medskip

    Next, we show that if $\pi$ contains $\pi_i\pi_j\pi_k\simeq 231$, its preimage $\tau$ contains either $1324$ or $\mu_{2413}$. Consider the following four cases for the ordering of $\pi_i$, $\pi_j$, and $\pi_k$ in $\tau$:
    \begin{itemize}
        \item Suppose $\pi_i$ occurs before $\pi_j$ in $\tau$, so that $\pi_i$ must be popped from the stack before $\pi_j$ enters. Say $\pi_i$ was popped when $\tau_p$ was being added to the stack because $\tau_p\pi_i\tau_r\simeq \un{23}1$. Then $\tau$ contains $\tau_r\pi_i\tau_p\pi_j\simeq 1324$. 
        \item Next, suppose $\tau$ contains $\pi_k\pi_j\pi_i$, so that $\pi_j$ must not pop before $\pi_i$ enters. Then, there must be some element less than $\pi_k$ between $\pi_i$ and $\pi_j$ in the stack. Choose $l$ to be the least index of such an element $\tau_l$. Henceforth, we consider the elements of $\tau$ with respect to the plot of $\pi_k\pi_j\tau_l\pi_i\simeq 2413$. By our choice of $\tau_l$, the squares $(2,0)$ and $(2,1)$ are empty. If there existed some element $\tau_s$ in the $(1,0)$ square, then before $\tau_l$ is added to the stack, all the elements between $\pi_j$ and $\tau_l$ as well as $\pi_j$ itself will be popped, a contradiction. Next, if the $(2,2)$ square is not empty, choose the least index $t$ corresponding to some $\tau_t$ in the $(2,2)$ square. Then before $\tau_t$ is added to the stack, all the elements between $\pi_j$ and $\tau_l$, as well as $\pi_j$ itself, will be popped, again a contradiction. Thus, $\tau$ contains $\pi_k\pi_j\tau_l\pi_i\simeq\mu_{2413}$. 
        \item Say $\tau$ contains $\pi_j\pi_i\pi_k$. Similarly to the previous case, $\pi_j$ must not pop before $\pi_i$ enters. However, both elements must pop before $\pi_k$ enters. Therefore, there exists an element $\pi_k'<\pi_i$ below $\pi_j$ as well as an element $\tau_l<\pi_k'$ between $\pi_i$ and $\pi_j$ in the stack. Repeat the previous argument, we have that $\tau$ contains $\pi_k'\pi_j\tau_l\pi_i\simeq\mu_{2413}$.
        \item Lastly, if $\tau$ contains $\pi_j\pi_k\pi_i$, the element $\pi_j$ must pop before $\pi_k$ enters but after $\pi_i$ enters, which is impossible. 
    \end{itemize}
\end{proof}

\medskip

\begin{proposition}
    It does not hold that $\sort(\sc_{\un{23}1})$ is a permutation class since $25314$ contains $2413$ but $\sc_{\un{23}1}(25314)=54132$ avoids $231$ while $\sc_{\un{23}1}(2413)=3142$ does not avoid $231$. 
\end{proposition}

\medskip

The first ten terms of $(|\sort_n(\sc_{\un{23}1})|)_{n\ge 0}$ are \[1, 1, 2, 6, 22, 90, 394, 1806, 8558, 41586.\] 
These terms match the OEIS sequence A006318, the Schr\"oeder numbers. We leave this enumeration as a conjecture. 

\medskip

\begin{conjecture} \label{scs231enumeration}
    The sorting class of $\sc_{\un{23}1}$ is enumerated by \[|\sort_n(\sc_{\un{23}1})|=S_{n-1}.\]
\end{conjecture}

\medskip

\subsection{The Sorting Classes of $\sc_{1\un{32}}$ and $\sc_{\un{12}3}$} \label{decomposition}

The classes $\sort_n(\sc_{1\un{32}})$ and $\sort_n(\sc_{\un{12}3})$ cannot be described succinctly with (mesh) pattern avoidance. We begin by showing the following property of $\sc_{1\un{32}}$. 

\medskip

\begin{lemma} \label{ssc132decreasing}
    For $\tau\in \sy_n$, its image $\pi$ under $\sc_{1\un{32}}$ is a decreasing sequence if and only if $\tau$ is increasing. In other words, $\pi=n(n-1)\cdots 1$ if and only if $\tau=12\cdots n$. 
\end{lemma}
\begin{proof}
    Since $12\cdots n$ avoids $\rev(1\un{32})=\un{23}1$, its image under $\sc_{1\un{32}}$ indeed is $n(n-1)\cdots 1$. To show that $\tau=12\cdots n$ is unique, suppose that $\pi$ is decreasing. As $\tau_1$ cannot be popped, we have $\tau_1=\pi_n=1$. Then, $\tau_2$ also cannot be popped, as it sits on top of $\tau_1=1$ in the stack. Thus, $\tau_2=\pi_{n-1}=2$. Repeating this argument shows the desired. 
\end{proof}

\medskip

Split $\tau$ into subsequences by $\tau=m_1a_1m_2a_2\cdots m_ka_k$, where $m_i$ is the $i$-th left-to-right minimum of $\tau$. Note that $m_k=1$. With \Cref{ssc132decreasing}, we characterize $\sort_n(\sc_{1\un{32}})$ in \Cref{ssc132characterization}.

\medskip

\begin{theorem} \label{ssc132characterization}
    The sorting class of $\sc_{1\un{32}}$ is composed of the permutations $\tau\in \sy_n$ such that $\rev(a_1)\rev(a_2)\cdots\rev(a_k)$ is decreasing.     
\end{theorem}
\begin{proof}
    We first show that none of the left-to-right minimums of $\tau$ are premature. Indeed, if $m_i$ is premature, then the term $\tau_j$ immediately below $m_i$ in the stack must be less than $m_i$. But $\tau_j$ comes before $m_i$ in $\tau$, contradicting the left-to-right minimality of $m_i$. 

    \medskip

    Let $b_i=\sc_{1\un{32}}(a_i)$. Since each term of $a_i$ is greater than $m_i$, the subsequence $a_i$ is sorted independently of the elements below $m_i$ in the stack. When $m_{i+1}$ is added, the remaining elements of $a_i$ still in the stack are all popped, as $m_{i+1}$ and $m_i$ form the $1\un{32}$ pattern with the element immediately above $m_i$ in the stack. It follows that \[\pi=\sc_{1\un{32}}(\tau)=b_1b_2\cdots b_km_km_{k-1}\cdots m_1.\]
    For $\pi$ to avoid $231$, the subsequence $b_1b_2\cdots b_k$ must be decreasing, by letting $m_k=1$ correspond to the $1$ in the pattern $231$. As a consequence of \Cref{ssc132decreasing}, we have that each $a_i$ is increasing. In particular, $b_i=\rev(a_i)$, so that $\rev(a_1)\rev(a_2)\cdots\rev(a_k)$ must be decreasing, as claimed. 

    \medskip

    The condition is sufficient, as $\pi$ then consists of a decreasing sequence followed by an increasing sequence, which clearly avoids $231$. 
\end{proof}

\medskip

\begin{proposition}
    It does not hold that $\sort(\sc_{1\un{32}})$ is a permutation class since $2413$ contains $132$ but $\sc_{1\un{32}}(2413)=4312$ avoids $231$ while $\sc_{1\un{32}}(132)=231$ does not avoid $231$. 
\end{proposition}

\medskip

\citet{Defant} studied the sorting class of $\sc_{\un{132}}$. Their conclusion is listed in \Cref{sc132machine}. 

\medskip

\begin{theorem}[{\citet[Proposition 4.1]{Defant}}] \label{sc132machine} 
    A permutation $\tau\in \sy_n$ is in $\sort_n(\sc_{\un{132}})$ if and only if $\rev(a_1)\rev(a_2)\cdots\rev(a_k)$ is decreasing. 
\end{theorem}

\medskip

It follows from \Cref{ssc132characterization} and \Cref{sc132machine} that $\sort_n(\sc_{1\un{32}})=\sort_n(\sc_{\un{132}})$. For $\tau\in \sort_n(\sc_{1\un{32}})$, we additionally note that $\sc_{1\un{32}}(\tau)=\sc_{\un{132}}(\tau)$, from the proof of \Cref{sc132machine} in \citet{Defant}. However, the two maps are distinct. For example, $\sc_{1\un{32}}(2431)=3412\ne 1342=\sc_{\un{132}}(2431)$. 

\medskip

\begin{theorem}[{\citet[Theorem 4.4]{Defant}}] \label{dyck132bijection} 
    The sorting class of $\sc_{\un{132}}$ is enumerated by \[|\sort_n(\sc_{\un{132}})|=C_n.\]
\end{theorem}

\medskip

\Cref{ssc132enumeration} follows from \Cref{dyck132bijection}, completing our investigation on $\sort(\sc_{1\un{32}})$. 

\medskip

\begin{corollary} \label{ssc132enumeration}
    The sorting class of $\sc_{1\un{32}}$ is enumerated by \[|\sort_n(\sc_{1\un{32}})|=C_n.\]
\end{corollary}

\medskip

Next, we describe $\sort_n(\sc_{\un{12}3})$, beginning with \Cref{un123lemma} and \Cref{s12av231}. 

\medskip

\begin{lemma}
    Consider $\tau\in \sy_n$ that is being processed by $\sc_{\un{12}3}$, and suppose $\tau_i$ is an element currently in the stack. Subsequently, after each time an element enters the stack, the stack will contain some element of value at most $\tau_i$. 
    \label{un123lemma}
\end{lemma}
\begin{proof}
    When $\tau_i$ is popped from the stack, the next element $\tau_p$ to be added to the stack must be less than $\tau_i$. 
\end{proof}

\medskip

\begin{lemma}[{\citet[p.6]{Defant}}] \label{s12av231} 
    The sorting class of $\sc_{12}$ is characterized by \[\sort_n(\sc_{12})=\av_n(213).\] 
\end{lemma}

\medskip

Now, let $\tau^a$ and $\tau^b$ be subsequences of $\tau$ so that $\tau=\tau^an\tau^b$. We use \Cref{un123lemma} and \Cref{s12av231} to characterize $\sort_n(\sc_{\un{12}3})$. 

\medskip

\begin{theorem}
    The sorting class of $\sc_{\un{12}3}$ consists of the permutations $\tau$ such that
    \begin{enumerate}
        \item each element in $\tau^a$ is greater than each element in $\tau^b$,
        \item $\tau^a\in\av(3\un{21},132)$, and
        \item $\tau^b\in\av(213)$. 
    \end{enumerate}
    \label{scs123characterization}
\end{theorem}
\begin{proof}
    Let $\pi=\sc_{\un{12}3}(\tau)$. Since $n$ cannot be popped from the stack, once it is added, the rest of the input sequence is essentially sorted by a $12$-avoiding stack. Thus, define $\pi^b=\sc_{12}(\tau^b)$. Let $\pi^a$ consist of the premature elements of $\sc_{\un{12}3}(\tau^a)$, and let $\pi^c$ be the remaining elements, so that $\sc_{\un{12}3}(\tau^a)=\pi^a\pi^c$. By following the definition of the stack-sorting map, we see that $\pi=\pi^a\pi^bn\pi^c$. 

    \medskip

    We first show that the listed conditions are necessary for $\pi$ to avoid $231$. Letting $n$ correspond to the $3$ in the pattern $231$, we see that each element of $\pi$ before $n$ must be less than each element of $\pi$ after $n$. In particular, each element of $\pi^a$ must be less than each element of $\pi^c$. \Cref{un123lemma} implies that the smallest element of $\tau^a$ cannot be premature, and thus will be in $\pi^c$. As a result, $\pi^a$ must be empty, from which condition (1) follows. Furthermore, $\tau^a$ must avoid $\av(3\un{21})$, so that $\pi^c=\rev(\tau^a)$. Since $\pi^c$ avoids $231$, we have that $\tau^a$ avoids $132$. Lastly, \Cref{s12av231} implies condition (3), since $\pi^b$ has to avoid $231$. 

    \medskip

    To check that the conditions are sufficient, suppose they all hold. Let $i$ be the index such that $\pi_i=n$. We show that each subsequence $\pi_j\pi_k\pi_l$ of $\pi$ avoids $231$. Clearly, we have that $j\ne i$ and $l\ne i$. If $j$ and $l$ are both greater than $i$ or both less than $i$, then conditions (2) and (3) ensure that $\pi_c$ and $\pi_b$ individually avoid $231$. Lastly, if $j<i<l$, it suffices to consider only the case $k=i$. Condition (1) finishes.
\end{proof}

\medskip

\begin{proposition}
    It does not hold that $\sort(\sc_{\un{12}3})$ is a permutation class since $4132$ contains $132$ but $\sc_{\un{12}3}(4132)=3214$ avoids $231$ while $\sc_{\un{12}3}(132)=231$ does not avoid $231$. 
\end{proposition}

\medskip

\Cref{3un21132enumeration} then allows us to enumerate $\sort_n(\sc_{\un{12}3})$ in \Cref{scs123enumeration}. 

\medskip

\begin{lemma} \label{3un21132enumeration}
    We have $|\av_n(3\un{21},132)|=2^{n-1}$ for $n\ge 1$. Additionally, $|\av_0(3\un{21},132)|=1$. 
\end{lemma}
\begin{proof}
    Consider $\tau\in \sy_n$ and let $\tau=\tau^an\tau^b$, as usual. We claim that $\tau\in\av_n(3\un{21},132)$ is equivalent to the following two conditions: 
    \begin{enumerate}
        \item $\tau^b=12\cdots (n-i)$, and
        \item $\tau^a\in \av(3\un{21},132)$. 
    \end{enumerate}

    \medskip
    
    Letting $n$ correspond to the $3$ in the pattern $3\un{21}$, we see that $\tau_b$ must be increasing. Furthermore, letting $n$ correspond to the $3$ in the pattern $132$, we have that each element of $\tau^b$ must be less than each term of $\tau^a$. Condition (1) follows. Condition (2) is clearly necessary. The conditions are sufficient, as $\tau^a$ and $\tau^b$ both individually avoid $3\un{21}$ and $132$. 

    \medskip

    Clearly, $|\av_0(3\un{21},132)|=|\av_1(3\un{21},132)|=1$. Now, let $\tau_i=n$. Due to condition (1), the index $i$ determines the values of the elements of $\tau^a$ and $\tau^b$. Those elements are then independently ordered by the two conditions. By summing over $i$, we have that \[|\av_n(3\un{21},132)|=\sum_{i=1}^{n}|\av_{n-i}(3\un{21},132)|.\]
    From this recurrence follows the lemma. 
\end{proof}

\medskip

Finally, we enumerate $\sort_n(\sc_{\un{12}3})$ with \Cref{scs123characterization} and \Cref{3un21132enumeration}.

\medskip

\begin{theorem} \label{scs123enumeration}
    The sorting class of $\sc_{\un{12}3}$ is enumerated by \[|\sort_n(\sc_{\un{12}3})|=C_{n-1}+\sum_{i=0}^{n-2}2^{n-2-i}C_i.\]
\end{theorem}
\begin{proof}
    Let $\tau=\tau^an\tau^b$, as usual. From \Cref{scs123characterization}, the index $i$ of $n$ determines $\tau^a$ and $\tau^b$. Summing over $i$, with \Cref{3un21132enumeration}, we have that 
    \begin{align*}    |\sort_n(\sc_{\un{12}3})|&=\sum_{i=1}^n|\av_{i-1}(3\un{21},132)|\cdot|\av_{n-i}(213)|\\ 
    &=C_{n-1} +\sum_{i=2}^n 2^{i-2} C_{n-i}=C_{n-1}+\sum_{i=0}^{n-2}2^{n-2-i}C_i,
    \end{align*}
    as desired. 
\end{proof}

\medskip

\subsection{The remaining maps $\sc_\s$} \label{leftoverpermutationclass}

In \Cref{avoidance,mesh,decomposition}, we fully characterized the sorting class of $\sc_\s$ for seven $\s\in \S_3$. In this section, we show that the five remaining sorting classes are not permutation classes. 

\medskip

\begin{proposition}
    For $\s\in \{\un{21}3,\un{31}2,2\un{13},2\un{31},3\un{12}$, it does not hold that $\sort(\sc_\s)$ is a permutation class. 
\end{proposition}
\begin{proof}
    We address each listed pattern in order.
    \begin{itemize}
        \item $4132$ contains $132$ but $\sc_{\un{21}3}(4132)=1234$ avoids $231$ while $\sc_{\un{21}3}(132)=231$ does not avoid $231$; 
        \item $3142$ contains $132$ but $\sc_{\un{31}2}(3142)=1243$ avoids $231$ while $\sc_{\un{31}2}(132)=231$ does not avoid $231$;
        \item $3142$ contains $132$ but $\sc_{2\un{13}}(3142)=4123$ avoids $231$ while $\sc_{2\un{13}}(132)=231$ does not avoid $231$;
        \item $361425$ contains $1324$ but $\sc_{2\un{31}}(361425)=165243$ avoids $231$ while $\sc_{2\un{31}}(1324)=3421$ does not avoid $231$; and 
        \item $3142$ contains $132$ but $\sc_{3\un{12}}(3142)=1243$ avoids $231$ while $\sc_{3\un{12}}(132)=231$ does not avoid $231$. 
    \end{itemize}
\end{proof}

\medskip

Recall that \citet[{Corollary 3.5}]{CCF} and \citet[{Theorem 2.2}]{Defant} addressed whether $\sort(\sc_\s)$ is a permutation class for classical and consecutive patterns $\s$, respectively. In summary, over $\s\in \S_3$, it holds that $\sort(\sc_\s)$ is a permutation class only for $\s\in\{321,1\un{23},3\un{21}\}$. 

\medskip 

\section{Dynamics of the maps $\sc_\s$} \label{dynamike}
In this section, we study some dynamical properties of $\sc_\s$. \Cref{preimages} is dedicated to finding the maximum number of preimages a permutation in $\sy_n$ can have under $\sc_\s$ for several $\s\in \S_3$. We additionally describe the unique permutation $\pi\in \sy_n$ fulfilling that maximum. \Cref{periodicpoints} deals with the periodic points of $\sc_\s$ for several $\s\in \S_3$. 

\medskip

\subsection{Maximum preimages of the maps $\sc_\s$} \label{preimages}

Define the \emph{movement sequence} of a permutation $\tau$ under $\sc_\s$ as a word composed of $n$ copies each of the letters $N$ and $X$ that represents the sequence of steps the stack of $\sc_\s$ takes as it processes $\tau$. We say such a movement sequence has size $n$. We use $N$ to denote a push operation (entry) and $X$ to denote a pop operation (exit), as in \citet{KatalinBerlow}. For example, the movement sequence of $\tau=514362$ under $\sc_{123}$ is $NNNXNNXXNXXX$ (see \Cref{exs123}). The following result from \citet{KatalinBerlow} bounds the number of preimages under $\sc_\s$. 

\medskip

\begin{theorem}[{\citet[Theorem 3.2]{KatalinBerlow}}] \label{maxpreimages} 
    Every permutation $\pi$ of length $n\ge 3$ has at most $C_{n-1}$ preimages under the map $\sc_\s$, for all patterns $\s\in \S_3$. That is, \[|\sc_\s^{-1}(\pi)|\le C_{n-1}.\]
\end{theorem}

\medskip

Fix $\pi$ and consider permutations $\tau$ satisfying $\sc_\s(\tau)=\pi$. Since $\tau_1$ cannot be premature, at any point (excepting the beginning and end) of the sorting process, there must have been at least one more push than pop to the stack of $\sc_\s$. Thus, each proper prefix of a movement sequence must have at least one more $N$ than $X$. \citet{KatalinBerlow} showed \Cref{maxpreimages} using this criteria on the possible movement sequences. Furthermore, should $|\sc_\s^{-1}(\pi)|=C_{n-1}$ for some $\pi\in \sy_n$ and $\s\in\S_3$, every valid movement sequence of size $n$ must correspond to an image $\tau$ of $\pi$. 

\medskip

Recall that by \Cref{complemma} we have $|\sc_\s^{-1}(\pi)|=|\sc_{\comp(\s)}^{-1}(\comp(\pi))|$. Therefore, it holds that $\max_{\pi\in \sy_n}|\sc_\s^{-1}(\pi)|=\max_{\pi\in \sy_n}|\sc_{\comp(\s)}^{-1}(\pi)|$. 

\medskip

\subsubsection{The maps $\sc_{\un{12}3}$ and $\sc_{\un{32}1}$}

As noted, it suffices to only consider the map $\sc_{\un{12}3}$. We show that the $C_{n-1}$ bound proposed in \Cref{maxpreimages} can be met. 

\medskip

\begin{lemma} \label{un123preimagesmeet}
    Let $\pi=(n-1)(n-2)\cdots 1n$. Then, $\sc_{\un{12}3}(\tau)=\pi$ if and only if $\tau_1=n$ and $\tau'\in\av_{n-1}(213)$, where $\tau'=\tau_2\tau_3\cdots\tau_n$. In particular, \[|\sc_{\un{12}3}^{-1}(\pi)|=|\av_{n-1}(213)|=C_{n-1}.\]
\end{lemma}
\begin{proof}
    As $\tau_1$ always stays in the stack, $\tau_1=n$ is necessary. Henceforth, $\tau'$ is essentially sorted by a stack avoiding $12$. The conclusion follows from \Cref{s12av231}. 
\end{proof}

\medskip

From \Cref{maxpreimages} and \Cref{un123preimagesmeet} it holds the following.

\medskip

\begin{corollary}
    We have
    \[\max_{\pi\in \sy_n}|\sc_{\un{12}3}^{-1}(\pi)|=\max_{\pi\in \sy_n}|\sc_{\un{32}1}^{-1}(\pi)|=C_{n-1}.\]
\end{corollary}

\medskip

We additionally show that the permutation in \Cref{un123preimagesmeet} is unique. 

\medskip

\begin{theorem}
    The maximum of $|\sc_{\un{12}3}^{-1}|$ is uniquely achieved by $\pi=(n-1)(n-2)\cdots 1n$. 
\end{theorem}
\begin{proof}
    Let $\pi$ satisfy $|\sc_{\un{12}3}^{-1}(\pi)|=C_{n-1}$. By the proof of \Cref{maxpreimages}, each valid movement sequence corresponds to some $\tau$ such that $\sc_{\un{12}3}^{-1}(\tau)=\pi$. Consider the movement sequence $N(NX)^{n-1}X$. Following the definition of the movement sequence, we have that $\tau_{i+1}\tau_i\tau_1\simeq123$ for each $i>1$. Thus $\tau_n<\tau_{n-1}<\cdots<\tau_1$, implying that $\tau=n(n-1)\cdots 21$. Finally, $\pi=\sc_{\un{12}3}(\tau)=(n-1)(n-2)\cdots 21n$, as desired. 
\end{proof}

\medskip

With \Cref{complemma}, we extend the above to $\sc_{\un{32}1}$. 

\medskip

\begin{corollary}
    The maximum of $|\sc_{\un{32}1}^{-1}|$ is uniquely achieved by $\pi=23\cdots n1$. 
\end{corollary}

\medskip

\subsubsection{The maps $\sc_{\un{23}1}$ and $\sc_{\un{21}3}$}

We primarily consider $\sc_{\un{23}1}$. We show that the $C_{n-1}$ bound proposed in \Cref{maxpreimages} can be met. 

\medskip

\begin{lemma} \label{un231preimagesmeet}
    Let $\pi=n(n-1)\cdots 1$. Then, $\sc_{\un{23}1}(\tau)=\pi$ if and only if $\tau_1=1$ and $\tau'\in\av_{n-1}(213)$, where $\tau'=\tau_2\tau_3\cdots\tau_n$. In particular, \[|\sc_{\un{23}1}^{-1}(\pi)|=|\av_{n-1}(213)|=C_{n-1}.\]
\end{lemma}
\begin{proof} 
    As $\tau_1$ always stays in the stack, $\tau_1=1$ is necessary. Henceforth, $\tau'$ is essentially sorted by a stack avoiding $12$. The conclusion follows from \Cref{s12av231}. 
\end{proof}

\medskip

From \Cref{maxpreimages} and \Cref{un231preimagesmeet} it holds the following.

\medskip

\begin{corollary}
    We have
    \[\max_{\pi\in \sy_n}|\sc_{\un{23}1}^{-1}(\pi)|=\max_{\pi\in \sy_n}|\sc_{\un{21}3}^{-1}(\pi)|=C_{n-1}.\]
\end{corollary}

\medskip

We additionally show that the permutation in \Cref{un231preimagesmeet} is unique. 

\medskip

\begin{theorem}
    The maximum of $|\sc_{\un{23}1}^{-1}|$ is uniquely achieved by $\pi=n(n-1)\cdots 1$. 
\end{theorem}
\begin{proof}
    Let $\pi$ satisfy $|\sc_{\un{23}1}^{-1}(\pi)|=C_{n-1}$. By the proof of \Cref{maxpreimages}, each valid movement sequence corresponds to some $\tau$ such that $\sc_{\un{23}1}^{-1}(\tau)=\pi$. Consider the movement sequence $N(NX)^{n-1}X$. Following the definition of the movement sequence, we have that $\tau_{i+1}\tau_i\tau_1\simeq231$ for each $i>1$. Thus $\tau_1<\tau_n<\tau_{n-1}<\cdots<\tau_2$, implying that $\tau=1n(n-1)\cdots 2$. Finally, $\pi=\sc_{\un{23}1}(\tau)=n(n-1)\cdots 1$, as desired. 
\end{proof}

\medskip

With \Cref{complemma}, we extend the above to $\sc_{\un{23}1}$. 

\medskip

\begin{corollary}
    The maximum of $|\sc_{\un{23}1}^{-1}|$ is uniquely achieved by $\pi=12\cdots n$. 
\end{corollary}

\medskip

\subsubsection{The maps $\sc_{2\un{31}}$ and $\sc_{2\un{13}}$}

We primarily consider $\sc_{2\un{31}}$, recalling that $\sc_{2\un{31}}=\sc_{231}$ by \Cref{free2stacks}. Consider the following special case of Theorem 3.7 from \citet{KatalinBerlow}. 

\medskip

\begin{theorem}[{\citet[Theorem 3.7]{KatalinBerlow}}] \label{classicalpreimagecap}
    Let $\s\in S_3$ be a classical pattern. There exists a permutation $\pi$ satisfying $|\sc_\s(\pi)|=C_{n-1}$ if and only if $\s_1$ and $\s_2$ are consecutive numbers. 
\end{theorem}

\medskip

Since $2$ and $3$ are consecutive, we have the following from \Cref{classicalpreimagecap}. 

\medskip

\begin{corollary} 
    We have
    \[\max_{\pi\in \sy_n}|\sc_{2\un{31}}^{-1}(\pi)|=\max_{\pi\in \sy_n}|\sc_{2\un{13}}^{-1}(\pi)|=C_{n-1}.\]
\end{corollary}

\medskip

We additionally show the following.

\medskip

\begin{theorem}
    The maximum of $|\sc_{2\un{31}}^{-1}|$ is uniquely achieved by $\pi=n(n-1)\cdots 1$. 
\end{theorem}
\begin{proof}
    Let $\pi$ satisfy $|\sc_{231}^{-1}(\pi)|=C_{n-1}$. By the proof of \Cref{maxpreimages}, each valid movement sequence corresponds to some $\tau$ such that $\sc_{231}^{-1}(\tau)=\pi$. Consider the movement sequence $N(NX)^{n-1}X$. Following the definition of the movement sequence, we have that $\tau_{i+1}\tau_i\tau_1\simeq231$ for each $i>1$. Thus $\tau_1<\tau_n<\tau_{n-1}<\cdots<\tau_2$, implying that $\tau=1n(n-1)\cdots 2$. Finally, $\pi=\sc_{231}(\tau)=n(n-1)\cdots 1$, as desired. 
\end{proof}

\medskip

With \Cref{complemma}, we extend the above to $\sc_{2\un{31}}$. 

\medskip

\begin{corollary}
    The maximum of $|\sc_{2\un{31}}^{-1}|$ is uniquely achieved by $\pi=12\cdots n$. 
\end{corollary}

\medskip

\subsubsection{The maps $\sc_{1\un{23}}$ and $\sc_{3\un{21}}$}

We primarily consider $\sc_{1\un{23}}$. 

\begin{conjecture} \label{1un23preimages}
    For $n\ge 2$, it holds that 
    \[\max_{\pi\in \sy_n}|\sc_{1\un{23}}^{-1}(\pi)|=\max_{\pi\in \sy_n}|\sc_{3\un{21}}^{-1}(\pi)|=2^{n-2}.\]
\end{conjecture}

\medskip

We first verify that a permutation achieving the aforementioned maximum exists. 

\medskip

\begin{lemma} \label{1un23achieves}
    Let $\pi=(n-1)(n-2)\cdots 1n$. Then, $\sc_{1\un{23}}(\tau)=\pi$ if and only if $\tau_1=n$ and $\tau'\in \av_{n-1}(132,213)$, where $\tau'=\tau_2\tau_3\cdots\tau_n$. In particular, \[|\sc_{1\un{23}}^{-1}(\pi)|=|\av_{n-1}(132,213)|=2^{n-2}.\]
\end{lemma}
\begin{proof}
    As the first element of $\tau$ always stays in the stack, $\tau_1=n$ is necessary. Henceforth, the entire stack (except $\tau_1$) will pop if there is ever added an element less than the second bottom-most element in the stack. Thus, $\tau'$ is the concatenation of increasing intervals of consecutive numbers, such that each interval begins with a left-to-right minimum of $\tau'$. In the notation of \citet{Defant}, $\tau'$ is \emph{reverse-layered}. It is known \cite{ReverseLayered} that the reverse-layered permutations precisely form the set $\av(132,213)$. We thus are able to enumerate $|\sc_{1\un{23}}^{-1}(\pi)|$ with \Cref{132213enumeration}. 
\end{proof}

\medskip

We partially resolve \Cref{1un23preimages} by induction on $n$, with the aim of showing that a permutation $\pi$ achieving the maximum value of $|\sc_{1\un{23}}^{-1}|$ is $\pi=(n-1)(n-2)\cdots 1n$. The bound claimed by \Cref{1un23preimages} would then follow by \Cref{1un23achieves}. The inductive base case $n=2$ is obvious. In what follows, for a permutation $\pi$, let $\pi'$ be obtained by swapping the elements $\pi_n$ and $\pi_n+1$ in $\pi$. The following two lemmas ultimately work to show \Cref{1un23lastelementsmoothed}. 

\medskip

\begin{lemma} \label{1un23smoothing}
    Consider a permutation $\pi\in \sy_n$ satisfying $\pi_n<n$ and $\pi_{n-1}\ne \pi_n+1$. Then, \[|\sc_{1\un{23}}^{-1}(\pi)|\le |\sc_{1\un{23}}^{-1}(\pi')|.\]
\end{lemma}
\begin{proof}
    Let $\tau$ be a preimage of $\pi$. Evidently, $\tau_1=\pi_n$. We claim that the elements $\pi_n$ and $\pi_n+1$ are never consecutive entries in the stack. Indeed, if they were, then $\pi_n+1$ would sit immediately above of $\pi_n$ in the stack, so that $\pi_n+1$ will never be popped from the stack. Therefore, $\pi_n+1$ and $\pi_n$ will be adjacent in $\pi$, contradicting our assumptions. Now, let $\tau'$ be obtained by swapping $\pi_n$ and $\pi_n+1$ in $\tau$. Note that $\pi_n$ and $\pi_n+1$ have the same size relative to the other elements of $\tau$, and that $\pi_n$ and $\pi_n+1$ cannot both be a part of the $1\un{23}$ pattern the stack checks for. It follows that the sequences of push and pop operations that sort $\tau$ and $\tau'$ are identical. In particular, $\sc_{1\un{23}}(\tau')=\pi'$ for all preimages $\tau$ of $\pi$. As the map $\tau\mapsto\tau'$ is injective, the lemma follows.  
\end{proof}

\medskip

\begin{lemma} \label{1un23lastadjacent}
    Consider $\pi\in \sy_n$ satisfying $\pi_{n-1}= \pi_n+1$. Then, \[|\sc_{1\un{23}}^{-1}(\pi)|\le 2^{n-3}.\] 
    In particular, \Cref{1un23achieves} implies that $\pi$ does not achieve the maximum value of $|\sc_{1\un{23}}^{-1}|$ over $\pi\in \sy_n$. 
\end{lemma}
\begin{proof}
    Let $\tau$ be a preimage of $\pi$. Evidently, $\tau_1=\pi_n$. We claim that $\tau_2=\pi_n+1$ if $\pi_{n-1}=\pi_n+1$. Note that $\pi_n+1$ must be immediately above $\pi_n$ in the stack at some point. Thus, when $\pi_n+1$ is the first element in the input sequence, all the terms in the stack (except $\pi_n$) must be popped. If $\tau_i$ is the element immediately above $\pi_n$ in the stack at that moment, we require $(\pi_n+1)\tau_i\pi_n\simeq 1\un{23}$, which is impossible. Thus, $\tau_2=\pi_n+1$. 

    \medskip

    In what follows, let $\tau^*\in S_{n-1}$ be obtained by removing $\pi_n+1$ from $\tau$ and shifting all the other elements of $\tau$ greater than $\pi_n$ down by $1$. Define $\pi^*\in S_{n-1}$ similarly. Note that there is a bijection between $\tau$ and $\tau^*$, as well as between $\pi$ and $\pi^*$. As $\pi_n+1$ will never be popped from the stack, we see that $\sc_{1\un{23}}(\tau^*)=\pi^*$. Therefore, \[|\sc_{1\un{23}}^{-1}(\pi)|=|\sc_{1\un{23}}^{-1}(\pi^*)|\le 2^{n-3},\] 
    by inductively using the statement of \Cref{1un23preimages}. 
\end{proof} 

\medskip

We are now ready to partially describe permutations that maximize $|\sc_{1\un{23}}^{-1}|$ by smoothing. 

\medskip

\begin{corollary} \label{1un23lastelementsmoothed}
    There exists $\pi$ achieving the maximum value of $|\sc_{1\un{23}}^{-1}|$ over $\sy_n$ such that the last element of $\pi$ is $\pi_n=n$. 
\end{corollary}
\begin{proof}
    Begin with some optimal $\pi\in \sy_n$. If $\pi_n<n$ and $\pi_{n-1}= \pi_n+1$, then $\pi$ is not optimal, by \Cref{1un23lastadjacent}. Thus $\pi_n<n$ and $\pi_{n-1}=\pi_n+1$, so that $\pi'$ has at least as many preimages as $\pi$, by \Cref{1un23smoothing}. Henceforth, we consider $\pi'$, and repeat these steps until we arrive at an optimal permutation such that its last element is $n$. 
\end{proof}

\medskip

We end by conjecturing that the permutation in \Cref{1un23achieves} is unique. 

\medskip

\begin{conjecture} \label{1un23uniquemaximum}
    The maximum of $|\sc_{1\un{23}}^{-1}|$ is uniquely achieved by $\pi=(n-1)(n-2)\cdots 1n$. Furthermore, the maximum of $|\sc_{3\un{21}}^{-1}|$ is uniquely achieved by $\pi=23\cdots n1$. 
\end{conjecture}

\medskip

\subsection{Periodic points of the maps $\sc_\s$} \label{periodicpoints}

Consider $\pi\in \sy_n$. In what follows, let $\pi^*\in S_{n-1}$ obtained by removing the element $1$ from $\pi$ and shifting the remaining elements of $\pi$ down by $1$. 

\medskip

\begin{theorem} \label{132231periodicpoints}
    For $\s\in \{\un{13}2,\un{23}1,1\un{32},2\un{31}\}$, the periodic points of $\sc_\s$ in $\sy_n$ are precisely $\av_n(132,231)$. For $n\ge 2$, these points have period $2$. 
\end{theorem}
\begin{proof}
    Let $\s=\un{23}1$, where a similar argument applies to the other patterns. We proceed by induction on $n$, where the $n=1$ case is trivial. Consider $\tau\in \sy_n$ and note that the elements $\tau_i=1$ and $\tau_j=2$ cannot be premature. Thus, if $\tau_j$ occurs before $\tau_i$ in $\tau$, their order will be flipped in $\sc_{\un{23}1}(\tau)$. Next, we suppose that $\tau_i$ occurs before $\tau_j$ in $\tau$. When $\tau_j$ enters the stack, all of the elements above $\tau_i$ in the stack will be popped. Therefore, $\tau_i$ and $\tau_j$ will be adjacent in the stack, and thus in $\pi=\sc_{\un{23}1}(\tau)$. But for permutations $\pi$ such that the element $1$ and $2$ are adjacent, we can easily check that $\sc_{\un{23}1}(\pi)^*=\sc_{\un{23}1}(\pi^*)$ and that $\pi^*\in\av_{n-1}(132,231)$ if and only if $\pi\in\av_n(132,231)$. Since $\pi^*$ is eventually sorted to $\av_{n-1}(132,231)$ by repeated applications of $\sc_{\un{23}1}$, $\pi$ will also eventually be sorted to $\av_n(132,231)$ by applications of $\sc_{\un{23}1}$. We conclude by induction. 

    \medskip

    Lastly, consider $\tau\in\av_n(132,231)$. Note that when $\tau$ is sorted by $\sc_{\un{23}1}$, the stack never triggers. Thus, $\sc_{\un{23}1}(\tau)=\rev(\tau)$, and indeed $\sc_{\un{23}1}^2(\tau)=\tau$. 
\end{proof}

\medskip

Furthermore, the proof of \Cref{132231periodicpoints} bounds the number of applications of $\sc_\s$ necessary for $\tau\in \sy_n$ to reach a periodic point. 

\medskip

\begin{corollary} \label{timetosortpart1}
    For $n\ge 2$, it holds that $\sc_\s^{2n-4}(\tau)\in\av_n(132,231)$ where $\tau\in \sy_n$ and $\s\in\{\un{13}2,\un{23}1,1\un{32},2\un{31}\}$. 
\end{corollary}
\begin{proof}
    We induct on $n$, with the $n=2$ case being apparent. From the proof of \Cref{132231periodicpoints}, the elements $1$ and $2$ are adjacent in $\sc_\s^2(\tau)$. It follows from $\sc_\s^{2n-6}(\tau^*)\in \av_{n-1}(132,231)$ that $\sc_\s^{2n-4}(\tau)\in \av_n(132,231)$. 
\end{proof}

\medskip

With \Cref{complemma}, we extend \Cref{132231periodicpoints} and \Cref{timetosortpart1} to the following. 

\medskip

\begin{corollary} \label{timetosortpart2}
    For $\s\in \{\un{31}2,\un{21}3,3\un{12},2\un{13}\}$, the periodic points of $\sc_\s$ of length $n$ are $\av_n(312,213)$. When $n\ge 2$, these points have period $2$. Furthermore, $\sc_\s^{2n-4}(\tau)\in\av_n(312,213)$ for each $\tau\in \sy_n$, for $n\ge 2$. 
\end{corollary}

\medskip

\section{Future Directions} \label{futuredirections}
In this work, we introduced the vincular-pattern-avoiding stack-sorting maps $\sc_\s$ as a generalization of \citet{CCF}'s and \citet{Defant}'s stack-sorting maps. A stack, that avoids $\s$, processes the input permutation $\tau$ in a right-greedy manner. For vincular patterns $\s$ of length $3$, we studied the sorting classes $\sort_n(\sc_\s)$, the preimages $\sc_\s^{-1}$, and the periodic points of $\sc_\s$. 

\medskip

In \Cref{sortingclass}, we identified and enumerated the sorting class $\sort_n(\sc_\s)$, the set of permutations in $\sy_n$ mapped onto the identity by $s\circ \sc_\s$ (see \Cref{sortingclassresults}). In \Cref{avoidance,mesh,decomposition}, we fully characterized the sorting class of $\sc_\s$ for seven $\s\in \S_3$. In \Cref{leftoverpermutationclass}, we discussed whether the five remaining sorting classes formed permutation classes. Here, we list computational results for $(|\sort_n(\sc_\s)|)_{n\ge 1}$ for those five sorting classes. 

\begin{footnotesize}
\begin{table}[h]
\centering
\begin{tabular}{c | c c c c c c c c c c c} 
    $\s\setminus n$ & 1 & 2 & 3 & 4 & 5 & 6 & 7 & 8 & 9 & OEIS & Permutation Class \\ 
    \hline\hline 
    \un{21}3 & 1 & 2 & 5 & 15 & 52 & 203 & 871 & 4017 & 19559 & Unknown & No \\
    \un{31}2 & 1 & 2 & 5 & 15 & 52 & 201 & 843 & 3764 & 17659 & \href{https://oeis.org/A202062}{A202062} & No \\
    2\un{13} & 1 & 2 & 5 & 16 & 62 & 273 & 1307 & 6626 & 35010 & Unknown & No \\
    2\un{31} & 1 & 2 & 6 & 23 & 102 & 496 & 2569 & 13934 & 78295 & Unknown & No \\
    3\un{12} & 1 & 2 & 5 & 15 & 52 & 201 & 843 & 3764 & 17659 & \href{https://oeis.org/A202062}{A202062} & No
\end{tabular}
\caption{Computational results for $(|\sort_n(\sc_\s)|)_{n\ge 1}$}
\label{sortingclassconjectures}
\end{table}
\end{footnotesize}

We note that, by \Cref{free2stacks}, we have $\sc_{\un{31}2}=\sc_{312}$, $\sc_{2\un{13}}=\sc_{213}$, and $\sc_{2\un{31}}=\sc_{231}$. The sorting classes of the maps $\sc_{312}$, $\sc_{213}$, and $\sc_{231}$ have been left unsolved by \citet{CCF} several years ago. We note that the classes $|\sort_n(\sc_{\un{31}2})|=|\sort_n(\sc_{3\un{12}})|$ are both enumerated by the OEIS sequence A202062. We thus list the following. 

\medskip

\begin{conjecture}
    The sorting classes of $\sc_{\un{31}2}$ and $\sc_{3\un{12}}$ are identical. That is, \[\sort_n(\sc_{312})=\sort_n(\sc_{\un{31}2})=\sort_n(\sc_{3\un{12}}).\] Furthermore, for $\tau\in \sort_n(\sc_{312})$, we have \[\sc_{312}(\tau)=\sc_{\un{31}2}(\tau)=\sc_{3\un{12}}(\tau).\]
\end{conjecture}

\medskip

Recall that we also have \Cref{scs231enumeration}, which states that $\sort_n(\sc_{\un{23}1})$ is enumerated by the Schr\"oder numbers (OEIS sequence A006318). 

\medskip

Next, we covered the maximum preimages of the maps $\sc_\s$ in \Cref{preimages}. We showed that, for several patterns $\s\in \S_3$, the maximum value of $|\sc_\s^{-1}|$ is $C_{n-1}$. We also stated \Cref{1un23preimages}, that $\max_{\pi\in \sy_n}|\sc_{1\un{23}}^{-1}(\pi)|=2^{n-2}$. We showed in \Cref{1un23lastelementsmoothed} that there exists a permutation $\pi\in \sy_n$ that maximizes $|\sc_{1\un{23}}(\pi)|$ such that $\pi_n=n$. We also conjectured that $\pi=(n-1)(n-2)\cdots 1n$ uniquely fulfills the maximum of $\sc_{1\un{23}}^{-1}(\pi)|$ in \Cref{1un23uniquemaximum}. Moreover, we conjecture the following with regards to the \emph{next-most} number of preimages of $\sc_{1\un{23}}$. 

\medskip

\begin{conjecture}
    The second-largest number of preimages under $\sc_{1\un{23}}$ that a permutation in $\sy_n$ can have is $2^{n-3}$, for $n\ge 3$. Furthermore, the number of permutations $\pi\in \sy_n$ satisfying $|\sc_{1\un{23}}^{-1}(\pi)|=2^{n-3}$ is $2n-2$. 
\end{conjecture}

\medskip

The maximum preimages of the maps $\sc_{\un{13}2}=\sc_{132}$ (see \Cref{free2stacks}) and $\sc_{1\un{32}}$ are unknown, though it is straightforward to show with \Cref{maxpreimages} that they are bounded strictly below $C_{n-1}$. We list computational results for $(\max_{\pi \in \sy_n}|\sc_\s^{-1}(\pi)|)_{n\ge 1}$ for $\s=\un{12}3$ and $\s=1\un{23}$ in the following table. 

\begin{footnotesize}
\begin{table}[h]
\centering
\begin{tabular}{c | c c c c c c c c c c c} 
    $\s\setminus n$ & 1 & 2 & 3 & 4 & 5 & 6 & 7 & OEIS \\ 
    \hline\hline 
    \un{13}2 & 1 & 1 & 2 & 3 & 7 & 16 & 37 & Unknown \\
    1\un{32} & 1 & 1 & 2 & 3 & 6 & 13 & 24 & Unknown \\
\end{tabular}
\caption{Computational results for $(|\sc_\s^{-1}|)_{n\ge 1}$}
\label{preimageconjectures}
\end{table}
\end{footnotesize}

We leave the preimages in \Cref{preimageconjectures} as directions for future study. 

\medskip

Lastly, we discussed the periodic points of the maps $\sc_\s$ in \Cref{periodicpoints}. We leave the identification of the periodic points of $\sc_\s$ for $\s\in \{\un{12}3,\un{32}1,1\un{23},3\un{21}\}$ open. Another area of future study is to evaluate the tightness of the bounds in \Cref{timetosortpart1} and \Cref{timetosortpart2}. \citet{Defant} showed, for each permutation $\tau\in \sy_n$, that $\sc_{132}^{n-1}(\tau)\in \av_n(132,231)$ and $\sc_{312}^{n-1}(\tau)\in \av_n(312,213)$. Furthermore, \citet{ChoiChoi} showed that $|\sc_{132}^{n-2}(\sy_n)|=|\sc_{312}^{n-2}(\sy_n)|=2^{n-1}+1$. Since it is known \cite{SimionSchmidt} that $|\av_n(132,231)|=|\av_n(312,213)|=2^{n-1}$, exactly one permutation in $\sy_n$ is not sorted by $n-2$ applications of $\sc_\s$ for $\s=132$ and $312$ to the respective periodic points. Of course, these results carry over to the maps $\sc_{\un{13}2}$ and $\sc_{\un{31}2}$, by \Cref{free2stacks}. Analogous results for the remaining maps $\sc_\s$ are left as future directions. 

\medskip

Another open direction is the study of the stack-sorting map $\sc_{\s_1,\s_2}$ that uses a stack avoiding both $\s_1$ and $\s_2$. When $\s_1$ and $\s_2$ are both classical patterns, the sorting class, preimages, and periodic points of $\sc_{\s_1,\s_2}$ have been evaluated \cite{TwoStacks, KatalinBerlow}. We jumpstart the diversification of $\s_1$ and $\s_2$ to include vincular patterns by listing some computational results for the enumeration of $\sort_n(\sc_{\s_1,\s_2})$. We include only the classes $\sort_n(\sc_{\s_1,\s_2})$ whose enumeration are known on the OEIS. As usual, $C_n$, $M_n$, and $S_n$ denote the $n$-th Catalan, Motzkin, and Schr\"oder number, respectively. 

\begin{footnotesize}
\begin{table}[h]
\centering 
\begin{tabular}{c c | c} 
    Enumerated by & OEIS & $(\s_1,\s_2)$\\ 
    \hline\hline 
    $2^n-n$ & \href{https://oeis.org/A000325}{A000325} & $(1\un{23},1\un{32})$ \\
    $C_n$ & \href{https://oeis.org/A000108}{A000108} & $(1\un{23},2\un{13})$, $(1\un{32},3\un{12})$, $(2\un{31},3\un{21})$, $(1\un{32},\un{31}2)$, \\
    & & $(3\un{12},\un{13}2)$, $(\un{12}3,\un{13}2)$, $(\un{13}2, \un{31}2)$, $(\un{23}1,\un{32}1)$ \\
    $M_n$ & \href{https://oeis.org/A001006}{A001006} & $(1\un{32},3\un{21})$, $(1\un{32},\un{32}1)$ \\
    $S_n$ & \href{https://oeis.org/A006318}{A006318} & $(1\un{32},2\un{31})$, $(1\un{32},\un{23}1)$, $(2\un{31},\un{13}2)$, $(\un{13}2,\un{23}1)$ \\
    $(n-1)2^{n-2}+1$ & \href{https://oeis.org/A005183}{A005183} & $(1\un{23},\un{13}2)$ \\
    $\binom{2n-2}{n-1}$ & \href{https://oeis.org/A000984}{A000984} & $(1\un{23},\un{23}1)$ \\
    $1+\sum_{i=1}^{n-1}(n-i)C_i$ & \href{https://oeis.org/A294790}{A294790} & $(1\un{32},\un{12}3)$ \\
    $\sum_{i=0}^{n-1}\binom{n-1}iC_i$ & \href{https://oeis.org/A007317}{A007317} & $(1\un{32},\un{13}2)$, $(3\un{12},\un{12}3)$, $(\un{12}3,\un{31}2)$ \\
    No closed form & \href{https://oeis.org/A202062}{A202062} & $(3\un{12},\un{31}2)$ \\
    $2^{n-1}$ & \href{https://oeis.org/A000079}{A000079} & $(\un{12}3,\un{32}1)$

\end{tabular}
\caption{Computational results for $|\sort_n(\sc_{\s_1,\s_2})|$}
\label{twostacksortingclassconjectures}
\end{table}
\end{footnotesize}

\medskip

We leave the proof of these results tabulated in \Cref{twostacksortingclassconjectures} as conjectures. 

\section*{Acknowledgements}
\noindent The author gratefully acknowledges Yunseo Choi for their guidance and suggestion of our direction of study. 

\bibliographystyle{abbrvnat}
\bibliography{references.bib}

\end{document}